# The new $(D^\alpha G)/G$ expansion method to solve Fractional KdV-Equations

**Uttam Ghosh**[1a] , **Susmita Sarkar**[1b] **and Shantanu Das** [2]
[1]Department of Applied Mathematics, University of Calcutta, Kolkata, India
[1a]email : uttam_math@yahoo.co.in
[1b]email : susmita62@yahoo.co.in
[2]Reactor Control Systems Design Section E & I Group BARC Mumbai India
email : shantanu@barc.gov.in

## Abstract

Fractional calculus of variation plays an important role to formulate the non-conservative physical problems. In this paper we use semi-inverse method and fractional variational principle to formulate the fractional order generalized Korteweg-deVries (KdV) equation with Jumarie type fractional derivative and proposed a new method to solve the non-linear fractional differential equation named as $D^\alpha G/G$ expansion method. Using this method we obtained the solutions of fractional order generalized KdV. The obtained solutions are more general compare to other method and the solutions are expressed in terms of the generalized hyperbolic, trigonometric functions and rational functions.

## 1.0 Introduction

Since the last decade the researchers of different fields of science (like in biology [4-5], physics [6-10], chaotic dynamics [11] control theory [12-14], material science [15], and in the physical systems that exhibit memory [1-3]); are using the fractional differential models. There is no suitable method using which one can easily establish the equations of motions of the problems in terms of the fractional differential equations, like we have for classical calculus case. Use of 'fractional variational principle' is almost logical method to establish the equation of motion of the system in terms of the fractional derivative. To formulate those problems, related to fractional dynamics first find the Lagrangian of the system, using the classical differential equations, or using semi-inverse method [16-18]. Then replace the classical order derivative by the fractional order derivative to from fractional Lagrangian of the system. Then fractional order variational principle method is used to establish the fractional order differential equations [19-21]. He's semi-inverse principle is the most easy method to construct the Lagrangian of the system from the given equation of motion [16-18].

The fractional calculus of variational problems was initially studied by different authors using Riemann-Liouville (RL) fractional derivative [22-24]. Malinowska *et al* [25] developed the fractional Euler-Lagrange's (EL) equations using the Caputo type fractional derivative. Since the Leibnitz rule is not valid for the fractional derivative of RL and Caputo types therefore there were some limitations in those formulations. In those cases both left and right RL and Caputo derivative is taken into consideration to represent the fractional EL equation. Jumarie [26] established the fractional EL equation in terms of modified RL derivative and formulated the fractional order Lagrangian and Hamiltonian; as in conjugation with classical calculus.

The equation of motion of non-conservative systems cannot be determined using the energy based approach [24, 25]. On the other hand Bauer [27] proved that it is impossible to use variational principle to derive a single linear dissipative equation of motion with constant coefficients. Thus the techniques of Lagrangian and Hamiltonian mechanics which are derived using variational principle are not valid. The derivative used in Bauer's theorem is classical (integer) order and thus the method fails. To overcome this problem Riewe [22-23] constructed Lagrangian using the non-integer order fractional derivative. Almost every process in nature is non-conservative; including the quantum processes [22] that can be expressed in terms of the non-integer order fractional derivative-or via use of fractional differential equations.

Due to the non-local nature of fractional calculus [19] the researchers of plasma physics gave attention to study electro-acoustic solitary waves (soliton) in plasma system with two different electron temperatures and with stationary ion, trapped electrons. Using the fractional differential models authors [9, 28-29] found the amplitude of soliton solution depends on the order of fractional order derivative.

On the other hand not only formulation of fractional differential equations have becomes an important mathematical-physics problem another current research is the solution of fractional differential equations (FDE). The FDE are the fractional generalization of the ordinary and partial differential equations. In past decades both mathematician and physicist made many contributions to solve the non-linear partial differential equations. Some of those methods are also valid for solving FDE only minor modifications are needed. Those methods are the Adomian Decomposition Method [12, 30, 31, 38], HPM [32-33], Differential Transform Method [34], Backlund transformation [35], Variational iteration method [36-37], and sub-equation methods [39]. The fractional sub equation method was proposed by Zhang *et al* [40] and is used by the authors [41-43] to find the exact analytical solutions of the non-linear fractional differential equations with modified RL derivative of Jumarie type [45].

Wang *et al* [44, 55] proposed $G'/G$ expansion method, is another important and effective method to find the travelling wave solution of non-linear partial differential equations, where the function $G$ satisfies the linear differential equation $G'' + \lambda G' + \mu G = 0$. In this paper we proposed the new method i.e. $(D^\alpha G)/G$ or $G^{(\alpha)}/G$ expansion method to construct the exact analytical solution of the Jumarie type non-linear fractional differential equations. Here the function $G$ satisfies the linear fractional differential equation that is following

$$D^{2\alpha}G + \lambda D^\alpha G + \mu G = 0 \tag{1.1}$$

where $D^{2\alpha} = D^\alpha D^\alpha \neq D^{\alpha+\alpha}$ is the Jumarie type differential operator. This method is new and not reported earlier anywhere. This new method is most effective to find more generalized solutions of the fractional order non-linear partial differential equation with Jumarie type derivative. The earlier method i.e. fractional sub-equation method uses the auxiliary function $\phi$ which satisfies fractional Riccati equation i.e. $D^\alpha \phi = \sigma + \phi^2$ that is non-linear. Whereas in this proposed new-method, the auxiliary function $G$ satisfies linear fractional differential equation (1.1). All the solutions obtained by fractional sub-equation method are covered by this new method, along with the new set of expressions; that we will elucidate in this paper.

Organization of the paper is as follows: In section-2 we describe the concept of fractional derivative and its development. Section-3 is devoted to formulate the fractional order generalized KdV equation. The new $(D^\alpha G)/G$ expansion method is described in section-4. The solution of fractional order KdV equation and modified KdV equation are found in section-5 and 6 respectively, applying this new approach. Finally conclusion is given in section-7; followed by list of references.

## 2.0 Basic concept of fractional order derivative and its development

There are different kind definitions of fractional derivative. The commonly used definitions are the Riemann-Liouville (RL) fractional derivative, Caputo derivative [1-2, 31, 46] and the Jumarie derivative [31, 45, 47]. The RL-definition is applicable for any functions which are one time integrable. Let $f(x)$ be any one time integrable function defined on $[a,b]$ then the 'left RL-definition' of fractional derivative of order $\alpha$ with $\alpha \in \mathbb{R}^+$ and $n \in \mathbb{N}$

$$_aD_x^\alpha f(x) = \frac{1}{\Gamma(n-\alpha)} \int_a^x f(t)(x-t)^{n-\alpha-1} dt, \quad n-1 \leq \alpha < n \tag{2.1}$$

The corresponding RL-fractional partial derivative is of the function $f(x,y)$ of several variable is defined [48-49] as

$$_aD_x^\alpha f(x,y) = \frac{1}{\Gamma(n-\alpha)} \int_a^x f(t,y)(x-t)^{n-\alpha-1} dt \qquad (2.2)$$

In 2006 Jumarie modified [45] the left RL definition of fractional derivative for continuous (but not necessarily differentiable) functions in the following form, assuming $f(0)$ as finite

$$\begin{aligned} f^{(\alpha)}(x) = D_x^\alpha f(x) &= \frac{1}{\Gamma(-\alpha)} \int_0^x (x-\xi)^{-\alpha-1} f(\xi) d\xi; & \alpha < 0 \\ &= \frac{1}{\Gamma(1-\alpha)} \frac{d}{dx} \int_0^x (x-\xi)^{-\alpha} \left( f(\xi) - f(0) \right) d\xi; & 0 < \alpha < 1 \qquad (2.3) \\ &= \left( f^{(\alpha-n)}(x) \right)^{(n)}; & n \leq \alpha < n+1 \quad n \geq 1 \end{aligned}$$

In terms of Jumarie fractional derivative the following rules are valid [31, 47, 49]

$$\begin{aligned} D_x^\alpha x^\gamma = \left( x^\gamma \right)^{(\alpha)} &= \frac{\Gamma(1+\gamma)}{\Gamma(1+\gamma-\alpha)} x^{\gamma-\alpha}; \quad \gamma > -1 \\ D_x^\alpha [f(x)g(x)] = \left( f(x)g(x) \right)^{(\alpha)} &= g(x) D_x^\alpha f(x) + f(x) D_x^\alpha g(x) \\ &= g(x) f^{(\alpha)}(x) + f(x) g^{(\alpha)}(x) \qquad (2.4) \\ D_x^\alpha \left[ f(g(x)) \right] = \left( f(g(x)) \right)^{(\alpha)} &= f_g'[g(x)] D_x^\alpha g(x) = D_g^\alpha f[g(x)] (g_x')^\alpha \\ &= f_g' \left( g(x) \right)^{(\alpha)} = \left( f(g(x)) \right)_g^{(\alpha)} \left( g_x' \right)^\alpha \end{aligned}$$

The fractional order Jumarie partial derivative of the function $f(x,y)$ is defined in the following way for the order $0 < \alpha < 1$, assuming $f(0,y)$ is finite function

$$D_x^\alpha f(x,y) = \frac{\partial^\alpha f(x,y)}{\partial x^\alpha} = \frac{1}{\Gamma(1-\alpha)} \frac{d}{dx} \int_0^x (x-t)^{-\alpha} \left( f(t,y) - f(0,y) \right) dt \qquad (2.5)$$

The Mittag-Leffler function [50] plays crucial role to solve the fractional order differential equations; as the exponential-functions does for the integer order classical differential equations. Jumarie [49] defined the fractional order (generalized or fractional) hyperbolic and trigonometric functions in terms of Mittag-Leffler functions in the following formulas

$$\tanh_\alpha(x^\alpha) = \frac{\sinh_\alpha(x^\alpha)}{\cosh_\alpha(x^\alpha)} \qquad \coth_\alpha(x^\alpha) = \frac{\cosh_\alpha(x^\alpha)}{\sinh_\alpha(x^\alpha)},$$

$$\sinh_\alpha(x^\alpha) = \frac{E_\alpha(x^\alpha) - E_\alpha(x^\alpha)}{2} \qquad \cosh_\alpha(x^\alpha) = \frac{E_\alpha(x^\alpha) + E_\alpha(x^\alpha)}{2} \qquad (2.6)$$

$$\tan_\alpha(x^\alpha) = \frac{\sin_\alpha(x^\alpha)}{\cos_\alpha(x^\alpha)} \qquad \cot_\alpha(x) = \frac{\cos_\alpha(x^\alpha)}{\sin_\alpha(x^\alpha)}$$

$$\sin_\alpha(x^\alpha) = \frac{E_\alpha(ix^\alpha) - E_\alpha(ix^\alpha)}{2i} \qquad \cos_\alpha(x^\alpha) = \frac{E_\alpha(ix^\alpha) + E_\alpha(ix^\alpha)}{2}$$

where $E_\alpha(z) = \sum_{k=0}^{\infty} \frac{z^\alpha}{\Gamma(1+k\alpha)}$ is the one parameter Mittag-Leffler function, defined in series form.

Again like the integer order calculus fractional order calculus play an important role to formulate the problems with non-conservative systems [25]. Riewe [22-23] was the pioneer in this field to develop the fractional calculus of variation. Agrawal [24] establish the fractional order Euler-Lagrange's equation in terms of the RL-fractional derivative in the form

$$\frac{\partial F}{\partial y(x)} + {}_xD_b^\alpha\left(\frac{\partial F}{\partial\left({}_aD_x^\alpha y(x)\right)}\right) + {}_aD_x^\beta\left(\frac{\partial F}{\partial\left({}_xD_b^\beta y(x)\right)}\right) = 0 \qquad (2.7)$$

where $F = F\left(x, y(x), {}_aD_x^\alpha y(x), {}_xD_b^\beta y(x)\right)$ is a continuous function. Similar equation was established using Caputo fractional derivative by Malinowska et al in [25]. Jumarie [26] establish the fractional Euler-Lagrange's equation for $F = F\left(x, y(x), {}_aD_x^\alpha y(x)\right)$ using modified fractional derivative in the form, that is following

$$\frac{\partial^\alpha F}{\partial y^\alpha} - \frac{d^\alpha}{dx^\alpha}\left(\frac{\partial F}{\partial\left({}_aD_x^\alpha y\right)}\right) = 0 \qquad (2.8)$$

In the next sections we will derive the Jumarie type fractional order generalized KdV equation using the above rule and will obtain the solutions of the equations in terms of generalized order hyperbolic, trigonometric functions and rational functions, with our new method.

## 3.0 Formulation of fractional order KdV equation
The generalized KdV can be taken in the form [51]
$$u_t + (au + hu^2 + cu^3 + f)u_x + bu_{xxx} = 0 \qquad (3.1)$$
Introducing the potential function $v$ defined as $u = v_x$ we get the potential equation of (3.1) in the following form

$$v_{xt} + (av_x + hv_x^2 + cv_x^3 + f)v_{xx} + bv_{xxxx} = 0 \qquad (3.2)$$

Using the variational principle [51] we consider the functional in the form

$$J(v) = \iint_{R \times T} v\left[c_1 v_{xt} + (ac_2 v_x + hc_3 v_x^2 + cc_6 v_x^3 + fc_4)v_{xx} + bc_5 v_{xxxx}\right] dx dt \qquad (3.3)$$

where $T$ stands for time unit and $R$ stands for space unit and $c_j$ for $j=1,2...,6$ are the Lagrange's Multipliers will be determine later. Using integration by parts with the conditions $v_x|_R = v_x|_T = v_t|_T = 0$ we obtain from (3.3) the following expression [51]

$$J(v) = \iint_{R \times T}\left[-c_1 v_t v_x - \tfrac{1}{2}ac_2 v_x^3 - \tfrac{1}{3}hc_3 v_x^4 - \tfrac{1}{4}cc_6 v_x^5 - fc_4 v_x^2 - \tfrac{1}{2}bc_5 v_{xx}^2\right] dx dt \qquad (3.4)$$

Thus Lagrangian of the equation (3.1) is given by following expression [51]

$$L(v_t, v_x, v_{xx}) = -c_1 v_t v_x - ac_2 \frac{v_x^3}{2} - hc_3 \frac{v_x^4}{3} - cc_6 \frac{v_x^5}{4} - fc_4 v_x^2 - bc_5 \frac{v_{xx}^2}{2}$$

The Euler-Lagrange's equation of the system is given by [51]

$$-\frac{\partial}{\partial t}\left(\frac{\partial L}{\partial v_t}\right) - \frac{\partial}{\partial x}\left(\frac{\partial L}{\partial v_x}\right) + \frac{\partial^2}{\partial x^2}\left(\frac{\partial L}{\partial v_{xx}}\right) = 0 \qquad (3.5)$$

Putting $L(v_t, v_x, v_{xx})$ from (3.4) in equation (3.5) we get the following

$$2c_1 v_{xt} + 3ac_2 v_x v_{xx} + 4hc_3 v_x^2 v_{xx} + 5cc_6 v_x^3 v_{xx} + 2fc_4 v_{xx} + bc_5 v_{xxxx} = 0 \qquad (3.6)$$

Comparing (3.6) and (3.1) we obtain $c_1 = \tfrac{1}{2}, c_2 = \tfrac{1}{3}, c_3 = \tfrac{1}{4}, c_4 = \tfrac{1}{2}, c_5 = 1, c_6 = \tfrac{1}{5}$. Thus we get the Lagrangian as expressed in following form [51]

$$L(v_t, v_x, v_{xx}) = -\frac{1}{2}v_t v_x - a\frac{v_x^3}{6} - h\frac{v_x^4}{12} - c\frac{v_x^5}{20} - f\frac{v_x^2}{2} - b\frac{v_{xx}^2}{2} \qquad (3.7)$$

Now we express the fractional counterpart of (3.7) considering the fractional Lagrangian with Jumarie type [26] fractional derivative, to write the following

$$L_\alpha(D_t^\alpha v, D_x^\alpha v, D_{xx}^{\alpha\alpha} v) = -\frac{1}{2}D_t^\alpha v D_x^\alpha v - a\frac{(D_x^\alpha v)^3}{6} - h\frac{(D_x^\alpha v)^4}{12} - c\frac{(D_x^\alpha v)^5}{20} - f\frac{(D_x^\alpha v)^2}{2} - b\frac{(D_{xx}^{\alpha\alpha} v)^2}{2} \qquad (3.8)$$

The functional of the space-time fractional differential equation is given by [26], i.e. in terms of fractional integration is

$$J_\alpha(v) = \iint_{R \times T} L(D_t^\alpha v, D_x^\alpha v, D_{xx}^{\alpha\alpha} v)(dx)^\alpha (dt)^\alpha \qquad (3.9)$$

Using the fractional variation principle we obtained fractional Euler-Lagrange's equation in the form

$$-D_t^\alpha\left(\frac{\partial L_\alpha}{\partial(D_t^\alpha[v_t])}\right) - D_x^\alpha\left(\frac{\partial L_\alpha}{\partial(D_x^\alpha[v_x])}\right) + D_{xx}^{\alpha\alpha}\left(\frac{\partial L_\alpha}{\partial(D_{xx}^{\alpha\alpha}[v_{xx}])}\right) = 0 \qquad (3.10)$$

That is also

$$-\frac{\partial^\alpha}{\partial t^\alpha}\left(\frac{\partial L_\alpha}{\partial (v_t)^{(\alpha)}}\right) - \frac{\partial^\alpha}{\partial x^\alpha}\left(\frac{\partial L_\alpha}{\partial (v_x)^{(\alpha)}}\right) + \frac{\partial^\alpha}{\partial x^\alpha}\frac{\partial^\alpha}{\partial x^\alpha}\left(\frac{\partial L_\alpha}{\partial (v_{xx})^{(\alpha\alpha)}}\right) = 0$$

After simplification we get the fractional order KdV equation that is

$$\tfrac{1}{2}\left(D_x^\alpha D_t^\alpha v + D_t^\alpha D_x^\alpha v\right) + \left(aD_x^\alpha v + h(D_x^\alpha v)^2 + c(D_x^\alpha v)^3 + f\right)D_{xx}^{\alpha\alpha}v + bD_{xxxx}^{\alpha\alpha\alpha\alpha}v = 0 \qquad (3.11)$$

In terms of Jumarie derivative [9, 52] the Schwartz relation is valid i.e. $D_x^\alpha D_t^\alpha v = D_t^\alpha D_x^\alpha v$. Now putting $u(x,y) = D_x^\alpha v$ equation (3.11) reduces to following

$$D_t^\alpha u + \left(au + hu^2 + cu^3 + f\right)D_x^\alpha u + bD_{xxx}^{\alpha\alpha\alpha}u = 0; \quad 0 < \alpha \leq 1. \qquad (3.12)$$

## 4.0 Description of new $(D^\alpha G)/G$ expansion method

In this section we shall describe the $(D^\alpha G)/G$ or $(G^{(\alpha)})/G$ method, where $G$ satisfies the linear $2\alpha$ – order linear fractional differential equation that is

$$D_\xi^{\alpha\alpha}G + \lambda D_\xi^\alpha G + \mu G = 0; \quad 0 < \alpha \leq 1 \qquad (4.1)$$

where $D_\xi^{\alpha\alpha}G = D_\xi^\alpha D_\xi^\alpha G$ and $\lambda, \mu$ are arbitrary constants and the partial derivatives are the Jumarie type fractional derivative [52]. Using the methodology describe in [53, 47] we find the solution of the fractional differential equation (4.1) in the following form

$$G(\xi) = \begin{cases} E_\alpha\left(-\tfrac{1}{2}\lambda\xi^\alpha\right)\left(A\cosh_\alpha\left(\tfrac{1}{2}\sqrt{\lambda^2-4\mu}\,\xi^\alpha\right) + B\sinh_\alpha\left(\tfrac{1}{2}\sqrt{\lambda^2-4\mu}\,\xi^\alpha\right)\right); & \lambda^2 - 4\mu > 0 \\ E_\alpha\left(-\tfrac{1}{2}\lambda\xi^\alpha\right)\left(A\cos_\alpha\left(\tfrac{1}{2}\sqrt{4\mu-\lambda^2}\,\xi^\alpha\right) + B\sin_\alpha\left(\tfrac{1}{2}\sqrt{4\mu-\lambda^2}\,\xi^\alpha\right)\right); & \lambda^2 - 4\mu < 0 \\ (A + B\xi^\alpha)E_\alpha\left(-\tfrac{1}{2}\lambda\xi^\alpha\right); & \lambda^2 - 4\mu = 0 \end{cases} \qquad (4.2)$$

where $A$ and $B$ are arbitrary constants; and $E_\alpha$ is one-parameter Mittag-Leffler function. Using (4.2) and the derivative formula defined in (2.4) following expressions are obtained

$$\frac{D^{\alpha}G}{G} = \begin{cases} -\dfrac{\lambda}{2} + \dfrac{\sqrt{\lambda^2-4\mu}}{2} \dfrac{\left(A\sinh_{\alpha}\left(\frac{1}{2}\sqrt{\lambda^2-4\mu}\,\xi^{\alpha}\right)()+B\cosh_{\alpha}\left(\frac{1}{2}\sqrt{\lambda^2-4\mu}\,\xi^{\alpha}\right)\right)}{\left(A\cosh_{\alpha}\left(\frac{1}{2}\sqrt{\lambda^2-4\mu}\,\xi^{\alpha}\right)()+B\sinh_{\alpha}\left(\frac{1}{2}\sqrt{\lambda^2-4\mu}\,\xi^{\alpha}\right)\right)} \\ \qquad\qquad\qquad\qquad\qquad\qquad\qquad\qquad :\lambda^2-4\mu>0 \\[6pt] -\dfrac{\lambda}{2} + \dfrac{\sqrt{4\mu-\lambda^2}}{2} \dfrac{\left(-A\sin_{\alpha}\left(\frac{1}{2}\sqrt{4\mu-\lambda^2}\,\xi^{\alpha}\right)+B\cos_{\alpha}\left(\frac{1}{2}\sqrt{4\mu-\lambda^2}\,\xi^{\alpha}\right)\right)}{\left(A\cos_{\alpha}\left(\frac{1}{2}\sqrt{4\mu-\lambda^2}\,\xi^{\alpha}\right)()+B\sin_{\alpha}\left(\frac{1}{2}\sqrt{4\mu-\lambda^2}\,\xi^{\alpha}\right)\right)} \\ \qquad\qquad\qquad\qquad\qquad\qquad\qquad\qquad :\lambda^2-4\mu<0 \\[6pt] -\dfrac{\lambda}{2} + \dfrac{B}{(A+B\xi^{\alpha})} \\ \qquad\qquad\qquad\qquad\qquad\qquad\qquad\qquad :\lambda^2-4\mu=0 \end{cases} \quad (4.3)$$

Consider the non-linear partial fractional differential equation of the following type

$$L(u, D_t^{\alpha}u, D_x^{\alpha}u, D_{xt}^{\alpha\alpha}u....) = 0 \qquad (4.4)$$

where $D_t^{\alpha}u = \frac{\partial^{\alpha}}{\partial t^{\alpha}}u$, $D_x^{\alpha}u = \frac{\partial^{\alpha}}{\partial x^{\alpha}}u$ and $u=u(x,t)$ is the unknown function. $L$ is a function of $u(x,t)$

**Step 1.** We use the following travelling wave transformation

$$\xi = x + ct \qquad (4.5)$$

and $u(\xi) = u(x,t)$. Then equation (4.4) reduce to the ordinary non-linear fractional differential equation form (using the same notation of fractional derivative)

$$L(u, D_{\xi}^{\alpha}u, D_{\xi}^{\alpha\alpha}u....) = 0 \qquad (4.6)$$

where $D_{\xi}^{\alpha}u = \frac{d^{\alpha}u}{d\xi^{\alpha}}$, $D_{\xi}^{\alpha\alpha}u = \frac{d^{\alpha}}{d\xi^{\alpha}}\left(\frac{d^{\alpha}u}{d\xi^{\alpha}}\right)$.

**Step 2.** The solution of (4.4) in the powers of $\frac{D^{\alpha}G}{G}$, we express in the following form

$$u(\xi) = \sum_{k=0}^{n} a_k \left(\frac{D^{\alpha}G}{G}\right)^k ; \quad a_n \neq 0$$

$$u(\xi) = a_0 + a_1\left(\frac{D^{\alpha}G}{G}\right) + a_2\left(\frac{D^{\alpha}G}{G}\right)^2 + .... + a_n\left(\frac{D^{\alpha}G}{G}\right)^n \qquad (4.7)$$

where $\frac{D^{\alpha}G}{G}$ is defined in (4.3) and $G(\xi)$ satisfies the equation (4.1), $a_0, a_1,...a_n$ are constants will be determined later.

**Step 3**. Use homogeneous balance principle to determine $n$. Then substitute (4.7) in (4.6) and compare the like power of $\frac{D^\alpha G}{G}$ to determine $a_0, a_1, ... a_n$. Finally put the values of $a_0, a_1, ... a_n$ in (4.7) to obtain the exact solutions.

## 5.0 Application of $(D^\alpha G)/G$ expansion method to fractional KdV equation

The space-time fractional KdV equation can be obtained in the form putting $a=1, h=c=f=0$ in (3.12), and we get the following (with writing $D_{xxx}^{\alpha\alpha\alpha}$ as $D_x^{\alpha\alpha\alpha}$)

$$D_t^\alpha u + u D_x^\alpha u + b D_x^{\alpha\alpha\alpha} u = 0 \tag{5.1}$$

Using travelling wave transformation defined in (4.5) i.e. $\xi = x + ct$, with $\xi'_t = c$ and $\xi'_x = 1$ and then using the (2.4) i.e. $D_x^\alpha\left[f(g(x))\right] = (f(g(x)))^{(\alpha)} = D_g^\alpha f[g(x)](g'_x)^\alpha = (f(g(x)))_g^{(\alpha)}(g'_x)^\alpha$; we get

$$D_t^\alpha u(x,t) = D_t^\alpha u(\xi) = \left(D_\xi^\alpha u(\xi)\right)(\xi'_t)^\alpha = c^\alpha D_\xi^\alpha u; \quad D_x^\alpha u(x,t) = D_x^\alpha u(\xi) = \left(D_\xi^\alpha u(\xi)\right)(\xi'_x)^\alpha = D_\xi^\alpha u$$

$$D_{xx}^{\alpha\alpha} u(x,t) = D_x^\alpha\left(D_x^\alpha u(\xi)\right) = D_x^\alpha\left(\left(D_\xi^\alpha u(\xi)\right)(\xi'_x)^\alpha\right) = D_x^\alpha\left(D_\xi^\alpha u\right) = \left(D_{\xi\xi}^{\alpha\alpha} u(\xi)\right)(\xi'_x)^\alpha = D_{\xi\xi}^{\alpha\alpha} u(\xi)$$

Doing the above steps once more we get $D_{xxx}^{\alpha\alpha\alpha} u(x,t) = D_{\xi\xi\xi}^{\alpha\alpha\alpha} u(\xi)$. Substituting these change of variables we get (5.1) as $c^\alpha D_\xi^\alpha u + u D_\xi^\alpha u + b D_\xi^{\alpha\alpha\alpha} u = 0$. Thus the equation (5.1) reduce to following equivalent expressions

$$c^\alpha D_\xi^\alpha u + u D_\xi^\alpha u + b D_\xi^{\alpha\alpha\alpha} u = 0; \quad D_\xi^\alpha\left[c^\alpha u + \tfrac{1}{2} u^2 + b D_\xi^{\alpha\alpha} u\right] = D_\xi^\alpha[K] \tag{5.2}$$

In (5.2) K is a constant, as, we have for Jumarie derivative of constant as zero i.e. $D_\xi^\alpha[K] = 0$. The term $\tfrac{1}{2} u^2$ appears (5.2); that get explained by using Jumarie formula of (2.4) i.e. $D_x^\alpha\left[f(g)\right] = (f(g))^{(\alpha)} = f'_g(g) D_x^\alpha[g]$. Applying this to $\tfrac{1}{2} u^2$ we get $D_\xi^\alpha\left(\tfrac{u^2}{2}\right) = \tfrac{1}{2}\left(D_\xi^\alpha u\right)\left(D_u^1 u^2\right) = u D_\xi^\alpha u$ i.e. by setting $f = u^2$, we have $f' = 2u$ and $g = u(\xi)$; we obtain the required result of (5.2). Here we assume $f = u^2$ as differentiable, which indeed is. Operating $D_\xi^{-\alpha}$ in the above equation of (5.2) we obtain the following expression

$$C_1 + c^\alpha u + \frac{u^2}{2} + b D_\xi^{\alpha\alpha} u = 0 \tag{5.3}$$

where $C_1$ is the arbitrary constant equals $-K$. Suppose the solution of the ordinary fractional differential equation of the form of (4.7) is $u(\xi) = \sum_{k=0}^{n} a_k \left(\frac{D^\alpha G}{G}\right)^k$, where $G(\xi)$ satisfies the equation (4.1). Using (4.7) i.e. $u(\xi) = a_0 + a_1 \left(\frac{D^\alpha G}{G}\right) + a_2 \left(\frac{D^\alpha G}{G}\right)^2 + .. + a_n \left(\frac{D^\alpha G}{G}\right)^n$ the following results is derived

$$(u(\xi))^2 = \left(\sum_{k=0}^{n} a_k \left(\frac{D^\alpha G}{G}\right)^k\right)^2 = \left(a_0 + a_1 \left(\frac{D^\alpha G}{G}\right) + a_2 \left(\frac{D^\alpha G}{G}\right)^2 + .. + a_n \left(\frac{D^\alpha G}{G}\right)^n\right)^2$$

$$u^2 = \sum_{k=0}^{n} a_n^2 \left(\frac{D^\alpha G}{G}\right)^{2n} + ...... \sum P_k$$

$$u^2 = a_0^2 + a_1^2 \left(\frac{D^\alpha G}{G}\right)^2 + a_2^2 \left(\frac{D^\alpha G}{G}\right)^4 + ...a_{n-1}^2 \left(\frac{D^\alpha G}{G}\right)^{2(n-1)} + a_n^2 \left(\frac{D^\alpha G}{G}\right)^{2n} + ... + \sum P_k$$

$$u^2 \sim a_n^2 \left(\frac{D^\alpha G}{G}\right)^{2n}$$

The $P_k$ being the product terms. Since to use homogeneous balance principle we need only the degree of highest order, i.e. $a_n^2 \left(\frac{D^\alpha G}{G}\right)^{2n}$ so other terms are not written explicitly, in last expression in above steps  We have the expression (4.1) i.e. $D_{\xi\xi}^{\alpha\alpha} G + \lambda D_\xi^\alpha G + \mu G = 0$ which gives $\left(D_{\xi\xi}^{\alpha\alpha} G\right)/G = -\lambda \left(D_\xi^\alpha G\right)/G - \mu = 0$. We do the following operation i.e. $D_\xi^\alpha \left(\frac{D_\xi^\alpha G}{G}\right)$, and write the following steps

$$D_\xi^\alpha \left(\frac{D_\xi^\alpha G(\xi)}{G(\xi)}\right) = \frac{G\left(D_\xi^\alpha \left[D_\xi^\alpha G\right]\right) - \left(D_\xi^\alpha G\right)^2}{G^2} = \frac{G\left(D_{\xi\xi}^{\alpha\alpha}\right) - \left(D_\xi^\alpha G\right)^2}{G^2}$$

$$= \frac{G\left(-\lambda D_\xi^\alpha G - \mu G\right) - \left(D_\xi^\alpha G\right)^2}{G^2}$$

$$= -\left(\frac{D_\xi^\alpha G}{G}\right)^2 - \lambda \left(\frac{D_\xi^\alpha G}{G}\right) - \mu$$

Using this above derived relation $D_\xi^\alpha \left(\frac{D_\xi^\alpha G}{G}\right) = -\left(\frac{D_\xi^\alpha G}{G}\right)^2 - \lambda \left(\frac{D_\xi^\alpha G}{G}\right) - \mu$ and the Jumarie fractional derivative formula (2.4) i.e. $D_x^\alpha \left[f(g(x))\right] = (f(g(x))^{(\alpha)} = f_g'[g(x)] D_x^\alpha g(x) = D_g^\alpha f[g(x)](g_x')^\alpha$ we have following derivation

$$D_\xi^\alpha\left(\left(\tfrac{D_\xi^\alpha G}{G}\right)^n\right) = \left(\left(\tfrac{D^\alpha G}{G}\right)^n\right)' \left(\tfrac{D_\xi^\alpha G}{G}\right)^{(\alpha)}$$

$$= \left(n\left(\tfrac{D_\xi^\alpha G}{G}\right)^{n-1}\right)\left(D_\xi^\alpha\left(\tfrac{D_\xi^\alpha G}{G}\right)\right)$$

$$= n\left(\tfrac{D_\xi^\alpha G}{G}\right)^{n-1}\left(-\left(\tfrac{D_\xi^\alpha G}{G}\right)^2 - \lambda\left(\tfrac{D_\xi^\alpha G}{G}\right) - \mu\right)$$

$$= n\left(-\left(\tfrac{D_\xi^\alpha G}{G}\right)^{n+1} - \lambda\left(\tfrac{D_\xi^\alpha G}{G}\right)^n - \mu\left(\tfrac{D_\xi^\alpha G}{G}\right)^{n-1}\right)$$

Implying the following approximate relations (with simplified notation of $D^\alpha$ instead $D_\xi^\alpha$ and taking the highest order of derivative)

$$D^\alpha u(\xi) = D^\alpha\left(a_0 + a_1\left(\tfrac{D^\alpha G}{G}\right) + a_2\left(\tfrac{D^\alpha G}{G}\right)^2 + \ldots + a_n\left(\tfrac{D^\alpha G}{G}\right)^n\right) \sim D^\alpha\left(a_n\left(\tfrac{D^\alpha G}{G}\right)^n\right)$$

$$\sim na_n\left(\tfrac{D^\alpha G}{G}\right)^{n+1}$$

Repeating above steps again, we get $D^\alpha D^\alpha u(\xi) \sim D^\alpha\left\{na_n\left(\tfrac{D^\alpha G}{G}\right)^{n+1}\right\} \sim n(n+1)a_n\left(\tfrac{D^\alpha G}{G}\right)^{n+2}$. Putting these above values in (5.3) i.e. $C_1 + c^\alpha u + \tfrac{1}{2}u^2 + bD_\xi^{\alpha\alpha}u = 0$ and comparing the highest power of $\tfrac{D^\alpha G}{G}$ from the highest order derivative term and the non-linear term we get $2n = n+2$., implying $n = 2$.

$$u(\xi) = a_0 + a_1\left(\frac{D^\alpha G}{G}\right) + a_2\left(\frac{D^\alpha G}{G}\right)^2 \tag{5.4}$$

From (5.4) we have the following

$$u^2 = a_0^2 + a_1^2\left(\frac{D^\alpha G}{G}\right)^2 + a_2^2\left(\frac{D^\alpha G}{G}\right)^4 + 2a_0a_1\left(\frac{D^\alpha G}{G}\right) + 2a_0a_2\left(\frac{D^\alpha G}{G}\right)^2 + 2a_1a_2\left(\frac{D^\alpha G}{G}\right)^3 \tag{5.5}$$

Doing the operation $D_\xi^\alpha$ on (5.4) and with using $D_\xi^\alpha\left(\tfrac{D_\xi^\alpha G}{G}\right) = -\left(\tfrac{D_\xi^\alpha G}{G}\right)^2 - \lambda\left(\tfrac{D_\xi^\alpha G}{G}\right) - \mu$, also using formula in (2.4) we write

$$D^\alpha u(\xi) = a_1\left(D^\alpha\left(\frac{D^\alpha G}{G}\right)\right) + a_2 D^\alpha\left(\left(\frac{D^\alpha G}{G}\right)^2\right)$$

$$= a_1\left(-\left(\tfrac{D^\alpha G}{G}\right)^2 - \lambda\left(\tfrac{D^\alpha G}{G}\right) - \mu\right) + a_2\left(2\left(\tfrac{D^\alpha G}{G}\right)\right)\left(D^\alpha\left(\tfrac{D^\alpha G}{G}\right)\right)$$

$$= a_1\left(-\left(\tfrac{D^\alpha G}{G}\right)^2 - \lambda\left(\tfrac{D^\alpha G}{G}\right) - \mu\right) + 2a_2\left(\frac{D^\alpha G}{G}\right)\left(-\left(\tfrac{D^\alpha G}{G}\right)^2 - \lambda\left(\tfrac{D^\alpha G}{G}\right) - \mu\right)$$

$$= -2a_2\left(\frac{D^\alpha G}{G}\right)^3 - (2a_2\lambda + a_1)\left(\frac{D^\alpha G}{G}\right)^2 - (2a_2\mu + a_1\lambda)\left(\frac{D^\alpha G}{G}\right) - a_1\lambda\mu$$

Doing one more fractional derivative $D^\alpha$ on above result we write (5.6), the following

$$D_\xi^{2\alpha} u(\xi) = D_\xi^\alpha D_\xi^\alpha u(\xi) = D^\alpha\left(-2a_2\left(\frac{D^\alpha G}{G}\right)^3 - (2a_2\lambda + a_1)\left(\frac{D^\alpha G}{G}\right)^2 - (2a_2\mu + a_1\lambda)\left(\frac{D^\alpha G}{G}\right) - a_1\lambda\mu\right)$$

$$\begin{aligned}
= &\, 6a_2\left(\frac{D^\alpha G}{G}\right)^4 + (10a_2\lambda + 2a_1)\left(\frac{D^\alpha G}{G}\right)^3 \\
&+ (8a_2\mu + 4a_2\lambda^2 + 3a_1\lambda)\left(\frac{D^\alpha G}{G}\right)^2 \\
&+ (6a_2\lambda\mu + 2a_1\mu + \lambda^2 a_1)\left(\frac{D^\alpha G}{G}\right) + \mu(2a_2\mu + \lambda a_1)
\end{aligned} \tag{5.6}$$

Putting (5.5) and (5.6) in (5.3), i.e. $C_1 + c^\alpha u + \tfrac{1}{2}u^2 + bD_\xi^{\alpha\alpha} u = 0$ and comparing the like powers of $\left(\tfrac{D^\alpha G}{G}\right)^m$ for $m = 0,1,2,3,4$; as we get up to $\left(\tfrac{D^\alpha G}{G}\right)^4$ from (5.6), we get following set of algebraic simultaneous equations

$$\left(\frac{D^\alpha G}{G}\right)^4 \quad : \frac{a_2^2}{2} + 6ba_2 = 0$$

$$\left(\frac{D^\alpha G}{G}\right)^3 \quad : a_1 a_2 + b(10\lambda a_2 + 2a_1) = 0$$

$$\left(\frac{D^\alpha G}{G}\right)^{2:} \quad : c^\alpha a_2 + \frac{1}{2}(a_1^2 + 2a_0 a_2) + b(8a_2\mu + 4a_2\lambda^2 + 3a_1\lambda) = 0 \qquad (5.7)$$

$$\left(\frac{D^\alpha G}{G}\right)^1 \quad : c^\alpha a_1 + a_0 a_1 + b(6a_2\lambda\mu + 2a_1\mu + \lambda^2 a_1) = 0$$

$$\left(\frac{D^\alpha G}{G}\right)^0 \quad : C_1 + c^\alpha a_1 + \frac{1}{2}a_0^2 + b\mu(2a_2\mu + \lambda a_1) = 0$$

Solving the above five equations we obtain the following values

$$\left.\begin{array}{l} a_1 = -12b\lambda, a_2 = -12b \\ c^\alpha = -a_0 - 8b\mu - b\lambda^2 \\ C_1 = -c^\alpha a_0 - \frac{1}{2}a_0^2 + 12b^2\mu(2\mu + \lambda) \end{array}\right\} \qquad (5.8)$$

and $a_0, \mu, \lambda$ are the arbitrary constants. Hence the solution (5.4) in terms of $\frac{D^\alpha G}{G}$ is written as

$$u(x,t) = u(\xi) = a_0 - 12b\lambda\left(\frac{D^\alpha G}{G}\right) - 12b\left(\frac{D^\alpha G}{G}\right)^2 \qquad (5.9)$$

Now considering the solutions as defined in (4.3), i.e. reproduced below

$$\frac{D^{\alpha}G}{G} = \begin{cases} -\dfrac{\lambda}{2} + \dfrac{\sqrt{\lambda^{2}-4\mu}}{2} \dfrac{\left(A\sinh_{\alpha}\left(\frac{1}{2}\sqrt{\lambda^{2}-4\mu}\,\xi^{\alpha}\right) + B\cosh_{\alpha}\left(\frac{1}{2}\sqrt{\lambda^{2}-4\mu}\,\xi^{\alpha}\right)\right)}{\left(A\cosh_{\alpha}\left(\frac{1}{2}\sqrt{\lambda^{2}-4\mu}\,\xi^{\alpha}\right) + B\sinh_{\alpha}\left(\frac{1}{2}\sqrt{\lambda^{2}-4\mu}\,\xi^{\alpha}\right)\right)} \\ \qquad\qquad\qquad\qquad\qquad\qquad\qquad\qquad\qquad : \lambda^{2} - 4\mu > 0 \\[1em] -\dfrac{\lambda}{2} + \dfrac{\sqrt{4\mu-\lambda^{2}}}{2} \dfrac{\left(-A\sin_{\alpha}\left(\frac{1}{2}\sqrt{4\mu-\lambda^{2}}\,\xi^{\alpha}\right) + B\cos_{\alpha}\left(\frac{1}{2}\sqrt{4\mu-\lambda^{2}}\,\xi^{\alpha}\right)\right)}{\left(A\cos_{\alpha}\left(\frac{1}{2}\sqrt{4\mu-\lambda^{2}}\,\xi^{\alpha}\right) + B\sin_{\alpha}\left(\frac{1}{2}\sqrt{4\mu-\lambda^{2}}\,\xi^{\alpha}\right)\right)} \\ \qquad\qquad\qquad\qquad\qquad\qquad\qquad\qquad\qquad : \lambda^{2} - 4\mu < 0 \\[1em] -\dfrac{\lambda}{2} + \dfrac{B}{(A+B\xi^{\alpha})} \\ \qquad\qquad\qquad\qquad\qquad\qquad\qquad\qquad\qquad : \lambda^{2} - 4\mu = 0 \end{cases}$$

The solution (5.9) i.e. $u(x,t) = u(\xi) = a_{0} - 12b\lambda\left(\frac{D^{\alpha}G}{G}\right) - 12b\left(\frac{D^{\alpha}G}{G}\right)^{2}$ can be written for different values of $\lambda$ and $\mu$. We will use some formulas derived in [56], where we established the following

$$\sin_{\alpha}(x+y)^{\alpha} = \sin_{\alpha}(x^{\alpha})\cos_{\alpha}(y^{\alpha}) + \cos_{\alpha}(x^{\alpha})\sin_{\alpha}(y^{\alpha})$$
$$\cos_{\alpha}(x+y)^{\alpha} = \cos_{\alpha}(x^{\alpha})\cos_{\alpha}(y^{\alpha}) - \sin_{\alpha}(x^{\alpha})\sin_{\alpha}(y^{\alpha})$$

Using the same procedure, as in [56] the similar type of formula for fractional hyperbolic functions are obtained as in following set

$$\sinh_{\alpha}(x+y)^{\alpha} = \sinh_{\alpha}(x^{\alpha})\cosh_{\alpha}(y^{\alpha}) + \cosh_{\alpha}(x^{\alpha})\sinh_{\alpha}(y^{\alpha})$$
$$\cosh_{\alpha}(x+y)^{\alpha} = \cosh_{\alpha}(x^{\alpha})\cosh_{\alpha}(y^{\alpha}) + \sinh_{\alpha}(x^{\alpha})\sinh_{\alpha}(y^{\alpha})$$

From the above four formulas the following expressions are deduced.

$$\tanh_{\alpha}(x+y)^{\alpha} = \frac{\sinh_{\alpha}(x+y)^{\alpha}}{\cosh_{\alpha}(x+y)^{\alpha}} = \frac{\sinh_{\alpha}(x^{\alpha})\cosh_{\alpha}(y^{\alpha}) + \cosh_{\alpha}(x^{\alpha})\sinh_{\alpha}(y^{\alpha})}{\cosh_{\alpha}(x^{\alpha})\cosh_{\alpha}(y^{\alpha}) + \sinh_{\alpha}(x^{\alpha})\sinh_{\alpha}(y^{\alpha})}$$

$$\tan_{\alpha}(x+y)^{\alpha} = \frac{\sin_{\alpha}(x+y)^{\alpha}}{\cos_{\alpha}(x+y)^{\alpha}} = \frac{\sin_{\alpha}(x^{\alpha})\cos_{\alpha}(y^{\alpha}) + \cos_{\alpha}(x^{\alpha})\sin_{\alpha}(y^{\alpha})}{\cos_{\alpha}(x^{\alpha})\cos_{\alpha}(y^{\alpha}) - \sin_{\alpha}(x^{\alpha})\sin_{\alpha}(y^{\alpha})}$$

The above expressions we will be using in subsequent cases and subsequent sections

**Case-I:** $\lambda^{2} - 4\mu > 0$

In this case the solution is

$$u(x,t) = a_0 - 12b\lambda \left( -\frac{\lambda}{2} + \frac{\sqrt{\lambda^2-4\mu}}{2} \frac{\left(A\sinh_\alpha\left(\frac{1}{2}\sqrt{\lambda^2-4\mu}\xi^\alpha\right) + B\cosh_\alpha\left(\frac{1}{2}\sqrt{\lambda^2-4\mu}\xi^\alpha\right)\right)}{\left(A\cosh_\alpha\left(\frac{1}{2}\sqrt{\lambda^2-4\mu}\xi^\alpha\right) + B\sinh_\alpha\left(\frac{1}{2}\sqrt{\lambda^2-4\mu}\xi^\alpha\right)\right)} \right)$$

$$-12b\left( -\frac{\lambda}{2} + \frac{\sqrt{\lambda^2-4\mu}}{2} \frac{\left(A\sinh_\alpha\left(\frac{1}{2}\sqrt{\lambda^2-4\mu}\xi^\alpha\right) + B\cosh_\alpha\left(\frac{1}{2}\sqrt{\lambda^2-4\mu}\xi^\alpha\right)\right)}{\left(A\cosh_\alpha\left(\frac{1}{2}\sqrt{\lambda^2-4\mu}\xi^\alpha\right) + B\sinh_\alpha\left(\frac{1}{2}\sqrt{\lambda^2-4\mu}\xi^\alpha\right)\right)} \right)^2$$

$$= a_0 + 3b\lambda^2 - 3b(\lambda^2-4\mu)\left( \frac{A\sinh_\alpha\left(\frac{1}{2}\sqrt{\lambda^2-4\mu}\xi^\alpha\right) + B\cosh_\alpha\left(\frac{1}{2}\sqrt{\lambda^2-4\mu}\xi^\alpha\right)}{A\cosh_\alpha\left(\frac{1}{2}\sqrt{\lambda^2-4\mu}\xi^\alpha\right) + B\sinh_\alpha\left(\frac{1}{2}\sqrt{\lambda^2-4\mu}\xi^\alpha\right)} \right)^2$$

(5.9)

$$u(x,t) = \begin{cases} a_0 + 3b\lambda^2 - 3b(\lambda^2-4\mu)\tanh_\alpha^2\left(\frac{1}{2}\sqrt{\lambda^2-4\mu}\xi^\alpha + \xi_0\right) \\ a_0 + 3b\lambda^2 - 3b(\lambda^2-4\mu)\coth_\alpha^2\left(\frac{1}{2}\sqrt{\lambda^2-4\mu}\xi^\alpha + \xi_1\right) \end{cases}$$

(5.10)

$$u(x,t) = \begin{cases} a_0 + 12b\mu + 3b(\lambda^2-4\mu)\text{sech}_\alpha^2\left(\frac{1}{2}\sqrt{\lambda^2-4\mu}\xi^\alpha + \xi_0\right) \\ a_0 + 12b\mu - 3b(\lambda^2-4\mu)\text{cosech}_\alpha^2\left(\frac{1}{2}\sqrt{\lambda^2-4\mu}\xi^\alpha + \xi_1\right) \end{cases}$$

(5.11)

where

$$\tanh_\alpha \xi_0 = \frac{B}{A}, \qquad \tanh_\alpha \xi_1 = \frac{A}{B} \qquad (5.12)$$

and

$$\xi = x + (-a_0 - 8b\mu - b\lambda^2)^{1/\alpha} t \qquad (5.13)$$

Thus the solution is expressed in terms of generalized hyperbolic functions.

**Case-II:** $\lambda^2 - 4\mu < 0$

In this case the solution is,

$$u(x,t) = a_0 - 12b\lambda\left(-\frac{\lambda}{2} + \frac{\sqrt{4\mu-\lambda^2}}{2}\frac{\left(-A\sin_\alpha\left(\frac{1}{2}\sqrt{4\mu-\lambda^2}\xi^\alpha\right) + B\cos_\alpha\left(\frac{1}{2}\sqrt{4\mu-\lambda^2}\xi^\alpha\right)\right)}{\left(A\cos_\alpha\left(\frac{1}{2}\sqrt{4\mu-\lambda^2}\xi^\alpha\right) + B\sin_\alpha\left(\frac{1}{2}\sqrt{4\mu-\lambda^2}\xi^\alpha\right)\right)}\right)$$

$$-12b\left(-\frac{\lambda}{2} + \frac{\sqrt{4\mu-\lambda^2}}{2}\frac{\left(-A\sin_\alpha\left(\frac{1}{2}\sqrt{4\mu-\lambda^2}\xi^\alpha\right) + B\cos_\alpha\left(\frac{1}{2}\sqrt{4\mu-\lambda^2}\xi^\alpha\right)\right)}{\left(A\cos_\alpha\left(\frac{1}{2}\sqrt{4\mu-\lambda^2}\xi^\alpha\right) + B\sin_\alpha\left(\frac{1}{2}\sqrt{4\mu-\lambda^2}\xi^\alpha\right)\right)}\right)^2$$

$$= a_0 + 3b\lambda^2 - 3b(4\mu-\lambda^2)\left(\frac{-A\sin_\alpha\left(\frac{1}{2}\sqrt{4\mu-\lambda^2}\xi^\alpha\right) + B\cos_\alpha\left(\frac{1}{2}\sqrt{4\mu-\lambda^2}\xi^\alpha\right)}{A\cos_\alpha\left(\frac{1}{2}\sqrt{4\mu-\lambda^2}\xi^\alpha\right) + B\sin_\alpha\left(\frac{1}{2}\sqrt{4\mu-\lambda^2}\xi^\alpha\right)}\right)^2$$

(5.14)

$$u(x,t) = \begin{cases} a_0 + 3b\lambda^2 - 3b(4\mu-\lambda^2)\tan_\alpha^2\left(\frac{1}{2}\sqrt{4\mu-\lambda^2}\xi^\alpha - \xi_0\right) \\ a_0 + 3b\lambda^2 - 3b(4\mu-\lambda^2)\cot_\alpha^2\left(\frac{1}{2}\sqrt{4\mu-\lambda^2}\xi^\alpha + \xi_1\right) \end{cases}$$

(5.15)

$$u(x,t) = \begin{cases} a_0 + 12b\mu - 3b(4\mu-\lambda^2)\sec_\alpha^2\left(\frac{1}{2}\sqrt{4\mu-\lambda^2}\xi^\alpha - \xi_0\right) \\ a_0 + 12b\mu - 3b(4\mu-\lambda^2)\csc_\alpha^2\left(\frac{1}{2}\sqrt{4\mu-\lambda^2}\xi^\alpha + \xi_1\right) \end{cases}$$

(5.16)

Where

$$\tan_\alpha \xi_0 = \frac{B}{A}; \quad \tan_\alpha \xi_1 = \frac{A}{B}; \quad \xi = x + (-a_0 - 8b\mu - b\lambda^2)^{1/\alpha} t$$

**Case-III:** $\lambda^2 - 4\mu = 0$

In this case the solution is

$$\therefore u(x,t) = a_0 - 12b\lambda\left(-\frac{\lambda}{2} + \frac{B}{(A+B\xi^\alpha)}\right) - 12b\left(-\frac{\lambda}{2} + \frac{B}{(A+B\xi^\alpha)}\right)^2$$

$$= a_0 + 3b\lambda^2 - \frac{12bB^2}{(A+B\xi^\alpha)^2} \qquad \xi = x + (-a_0 - 8b\mu - b\lambda^2)^{1/\alpha} t$$

(5.17)

The solutions are found in (5.9)-(5.11) and (5.14)-(5.17) are the generalized solution of fractional order KdV equation. The solutions (5.10), (5.15) and (5.17) are found when fractional sub-equation methods are used. Other solutions are new set appearing in this method. When $\alpha \to 1$ the solution reduce to the solution by Wang *et al* [55].

# 6.0 Application of $(D^\alpha G)/G$ expansion method to fractional modified KdV equation (mKdV-equation)

We will be using the same steps, as done for KdV equation earlier for mKdV equation in this section. Consider the mKdV equation in the form, that is obtained setting $a = c = f = 0$ and $h = 1$ in (3.12)

$$D_t^\alpha u - u^2 D_x^\alpha u + b D_x^{\alpha\alpha\alpha} u = 0 \tag{6.1}$$

Using travelling wave transformation defined in (4.5) the equation (6.1) reduces to following (as we demonstrated for KdV equation in previous section)

$$c^\alpha D_\xi^\alpha u - u^2 D_\xi^\alpha u + b D_\xi^{\alpha\alpha\alpha} u = 0$$

Using the derivative formula defined in (2.4), and with the steps we did for (5.2), we write the following form for the above equation

$$D_\xi^\alpha \left( c^\alpha u - \frac{u^3}{3} + b D_\xi^{\alpha\alpha} u \right) = D_\xi^\alpha [K] \tag{6.2}$$

Operating $D_\xi^{-\alpha}$ in both sides of (6.2) the above can be written in the following form,

$$c_1 + c^\alpha u - \frac{u^3}{3} + b D_\xi^{\alpha\alpha} u = 0 \tag{6.3}$$

where $c_1$ is the arbitrary constant will be determine later. Now consider the solution of (6.3) in the form defined in (4.7) as $u(\xi) = \sum_{k=0}^{n} a_k \left( \frac{D^\alpha G}{G} \right)^k$ and considering the homogeneous balance principle we obtain the following

$$u(\xi) = a_0 + a_1 \left( \frac{D^\alpha G}{G} \right) \tag{6.4}$$

where $G(\xi)$ satisfies the linear fractional differential equation (4.1) and $a_0, a_1$ are constant to be determine later. From (6.3) we have the following,

$$D_\xi^{\alpha\alpha} u = a_1 \left( 2\left(\frac{D^\alpha G}{G}\right)^3 + 3\lambda \left(\frac{D^\alpha G}{G}\right)^2 + (\lambda^2 + 2\mu)\left(\frac{D^\alpha G}{G}\right) + \lambda\mu \right) \tag{6.5}$$

and

$$u^3 = a_0^3 + a_1^3 \left(\frac{D^\alpha G}{G}\right)^3 + 3a_0^2 a_1 \left(\frac{D^\alpha G}{G}\right) + 3a_0 a_1^2 \left(\frac{D^\alpha G}{G}\right)^2 \tag{6.6}$$

Putting (6.4-6.6) in (6.3) i.e. $c_1 + c^\alpha u - \frac{1}{3}u^3 + bD_\xi^{\alpha\alpha} u = 0$ and comparing the comparing the coefficients of $\left(\frac{D^\alpha G}{G}\right)^m$ for $m = 0, 1, 2, 3$ (largest power is three in this case from (6.5) and (6.6)) we get the following set of equations

$$\begin{aligned}
\left(\frac{D^\alpha G}{G}\right)^3 & : -\frac{1}{3}a_1^3 + 2a_1 b = 0 \\
\left(\frac{D^\alpha G}{G}\right)^{2:} & : -a_1^2 a_0 + 3a_1 \lambda b = 0 \\
\left(\frac{D^\alpha G}{G}\right)^1 & : c^\alpha a_1 - a_0^2 a_1 + a_1(\lambda^2 + 2\mu)b = 0 \\
\left(\frac{D^\alpha G}{G}\right)^0 & : c_1 + c^\alpha a_0 - \frac{1}{3}a_0^3 + a_1 b\mu\lambda = 0
\end{aligned} \quad (6.7)$$

Solving the equations in (6.7) we obtain, the following values of the constants as

$$a_0 = \pm \frac{\lambda}{2}\sqrt{6b}, \quad a_1 = \pm\sqrt{6b}, \quad c^\alpha = \frac{1}{2}\lambda^\alpha b - 2b\mu, \quad c_1 = 0 \quad (6.8)$$

Hence the solution is given by $u(\xi) = \pm \frac{\lambda}{2}\sqrt{6b} \pm \sqrt{6b}\left(\frac{D^\alpha G}{G}\right)$. The solutions for different values of $\lambda$ and $\mu$ are given below.

**Case-I:** $\lambda^2 - 4\mu > 0$

In this case the solution is

$$\begin{aligned}
u(x,t) &= \pm \frac{\lambda}{2}\sqrt{6b} \\
&\pm \sqrt{6b}\left(-\frac{\lambda}{2} + \frac{\sqrt{\lambda^2 - 4\mu}}{2}\left(\frac{A\sinh_\alpha\left(\frac{1}{2}\sqrt{\lambda^2 - 4\mu}\xi^\alpha\right) + B\cosh_\alpha\left(\frac{1}{2}\sqrt{\lambda^2 - 4\mu}\xi^\alpha\right)}{A\cosh_\alpha\left(\frac{1}{2}\sqrt{\lambda^2 - 4\mu}\xi^\alpha\right) + B\sinh_\alpha\left(\frac{1}{2}\sqrt{\lambda^2 - 4\mu}\xi^\alpha\right)}\right)\right) \\
&= \pm\frac{1}{2}\sqrt{6b}\sqrt{\lambda^2 - 4\mu}\left(\frac{A\sinh_\alpha\left(\frac{1}{2}\sqrt{\lambda^2 - 4\mu}\xi^\alpha\right) + B\cosh_\alpha\left(\frac{1}{2}\sqrt{\lambda^2 - 4\mu}\xi^\alpha\right)}{A\cosh_\alpha\left(\frac{1}{2}\sqrt{\lambda^2 - 4\mu}\xi^\alpha\right) + B\sinh_\alpha\left(\frac{1}{2}\sqrt{\lambda^2 - 4\mu}\xi^\alpha\right)}\right)
\end{aligned} \quad (6.9)$$

$$u(x,t) = \begin{cases} \pm \dfrac{1}{2}\sqrt{6b}\sqrt{\lambda^2-4\mu}\ \tanh_\alpha\left(\tfrac{1}{2}\sqrt{\lambda^2-4\mu}\,\xi^\alpha + \xi_0\right) \\ \pm \dfrac{1}{2}\sqrt{6b}\sqrt{\lambda^2-4\mu}\ \coth_\alpha\left(\tfrac{1}{2}\sqrt{\lambda^2-4\mu}\,\xi^\alpha + \xi_1\right) \end{cases} \quad (6.10)$$

where $\tanh_\alpha \xi_0 = B/A$, $\tanh_\alpha \xi_1 = A/B$ and $\xi = x + \left(\tfrac{1}{2}\lambda^\alpha b - 2b\mu\right)^{1/\alpha} t$. Thus the solution is expressed in terms of generalized hyperbolic functions.

**Case-II:** $\lambda^2 - 4\mu < 0$

In this case the solution is,

$$u(x,t) = \pm \dfrac{\lambda}{2}\sqrt{6b}$$

$$\pm \sqrt{6b}\left(-\dfrac{\lambda}{2} + \dfrac{\sqrt{4\mu-\lambda^2}}{2}\dfrac{\left(-A\sin_\alpha\left(\tfrac{1}{2}\sqrt{4\mu-\lambda^2}\,\xi^\alpha\right)+B\cos_\alpha\left(\tfrac{1}{2}\sqrt{4\mu-\lambda^2}\,\xi^\alpha\right)\right)}{\left(A\cos_\alpha\left(\tfrac{1}{2}\sqrt{4\mu-\lambda^2}\,\xi^\alpha\right)+B\sin_\alpha\left(\tfrac{1}{2}\sqrt{4\mu-\lambda^2}\,\xi^\alpha\right)\right)}\right) \quad (6.11)$$

$$= \pm \dfrac{1}{2}\sqrt{6b}\sqrt{4\mu-\lambda^2}\ \dfrac{\left(-A\sin_\alpha\left(\tfrac{1}{2}\sqrt{4\mu-\lambda^2}\,\xi^\alpha\right)+B\cos_\alpha\left(\tfrac{1}{2}\sqrt{4\mu-\lambda^2}\,\xi^\alpha\right)\right)}{\left(A\cos_\alpha\left(\tfrac{1}{2}\sqrt{4\mu-\lambda^2}\,\xi^\alpha\right)+B\sin_\alpha\left(\tfrac{1}{2}\sqrt{4\mu-\lambda^2}\,\xi^\alpha\right)\right)}$$

$$u(x,t) = \begin{cases} \pm \dfrac{1}{2}\sqrt{6b}\sqrt{4\mu-\lambda^2}\ \tan_\alpha\left(\tfrac{1}{2}\sqrt{4\mu-\lambda^2}\,\xi^\alpha - \xi_0\right) \\ \pm \dfrac{1}{2}\sqrt{6b}\sqrt{4\mu-\lambda^2}\ \cot_\alpha\left(\tfrac{1}{2}\sqrt{4\mu-\lambda^2}\,\xi^\alpha + \xi_1\right) \end{cases} \quad (6.12)$$

where $\tan_\alpha \xi_0 = B/A$ and $\tan_\alpha \xi_1 = A/B$ and $\xi = x + \left(\tfrac{1}{2}\lambda^\alpha b - 2b\mu\right)^{1/\alpha} t$

**Case-III:** $\lambda^2 - 4\mu = 0$

In this case the solution is following, with $\omega$ as constant

$$\therefore u(x,t) = \pm \dfrac{\lambda}{2}\sqrt{6b} \pm \sqrt{6b}\left(-\dfrac{\lambda}{2} + \dfrac{B}{(A+B\xi^\alpha)}\right)$$

$$= \pm \dfrac{\sqrt{6b}}{\xi^\alpha + \omega}, \qquad \omega = \dfrac{A}{B} \quad (6.13)$$

The solutions are found in (6.9)-(6.10) and (6.11)-(6.13) are the generalized solution of space-time fractional order KdV equation. The solutions (6.10), (6.12) and (6.13) are found when fractional sub-equation methods are used. Other solutions are new set in this new method. When $\alpha \to 1$ the solutions reduce to the solution obtained by Wang *et al* [55].

## 7.0 Conclusions

In this paper we proposed a new method to solve the named as $(D^\alpha G)/G$ or $G^{(\alpha)}/G$ expansion method. Using this method we have found the exact analytical solution of fractional order KdV and the Modified KdV (mKdV) equation; which are fractional non-linear differential equations. The obtained solutions are expressed in terms of fractional order hyperbolic functions and fractional order trigonometric functions and the rational functions. Using this new proposed method we obtained more general solutions, compared to other method used earlier. Those obtained solutions reduce to the solutions obtained by $G'/G$ expansion method when $\alpha$ tends to one. This method is new and not reported earlier anywhere and this method is more effective to find more generalized solutions of the fractional order non-linear partial differential equation with Jumarie type derivative. In earlier method like-fractional sub-equation method the auxiliary function satisfies fractional Riccati equation that is non-linear, whereas in this method, the auxiliary function $G$ satisfies linear fractional differential equation. All the solutions obtained by earlier used fractional sub-equation method are covered by this method, along with the new set of expressions; that we presented in this paper.

## Reference

[1] I. Podlubny. Fractional Differential Equations. Mathematics in Science and Engineering, vol. 198. Academic Press, San Diego (1999)
[2] K. S. Miller, B. Ross. An Introduction to the Fractional Calculus and Fractional Differential Equations. Wiley, New York (1993)
[3] R. Hilfer. Applications of Fractional Calculus in Physics. World Scientific Publishing, River Edge (2000)
[4] M. Alquran, K. Al-Khaled, J. Chattapadhyay. Analytical solution of fractional population diffusion model: Residual power series. Non-Linear Studies. 22(1). 2015. 31-39.
[5] A. A. M. Arafa. Fractional differential equations in description of Bacterial growth. 21(3). 2013. Differ. Equ. Dyn. Syst.205-214.
[6] Parada F. J. V., Tapia J. A. O. and Ramirez J. A., 2007, Effective medium equations for fractional Fick's law in porous media, *Physica A*, 373, 339–353.
[7] Y. Chen, H. L. An. Numerical solution of Coupled Burger equations with time and space fractional derivative. Applied Mathematics Computation. 200. 2008. 87-95.
[8] X. Li and W. Chen. Analytical study of fractional anomulus diffusion in a half-plane. Journal of Physics A: Mathematical and Theoretical. 2010.43(49). 1-11.
[9] J. Banerjee, U. Ghosh , S. Sarkar and S. Das. A Study of Fractional Schrödinger Equation-composed via Jumarie fractional derivative. (in press Pramana - Journal of physics).

[10] Torvik P. J. and Bagley R. L., 1984, On the appearance of the fractional derivative in the behavior of real materials, *Transactions of the ASME*, 51, 294–298.

[11] Zaslavsky GM (2008) Hamiltonian chaos and fractional dynamics. Reprint of the 2005 original. Oxford University Press, Oxford

[12] S. Das. Functional Fractional Calculus 2$^{nd}$ Edition, Springer-Verlag 2011 and Functional Fractional Calculus for system identification & controls 1$^{st}$ Edition Springer-Verlag 2007.

[13] Axtell M. and Bise E. M., 1990, Fractional calculus applications in control systems, *Proc. of the IEEE Nat. Aerospace and Electronics Conf.*, New York, 563–566.

[14] Dorcak L. 1994, *Numerical Models for the Simulation of the Fractional-OrderControl Systems*, *UEF-04-94*, The Academy of Sciences, Inst. of Experimental Physic, Kosice, Slovakia.

[15] Mainardi F (2010) Fractional calculus and waves in linear viscoelasticity. Imperial College Press, London.

[16] J. He, Semi-inverse method of establishing generalized variational principles for fluid mechanics with emphasis on turbo-machinery aerodynamics, *Int. J. Turbo Jet-Engines* 14(1), 1997, pp. 23–28.

[17] J. He, Variational principles for some nonlinear partial differential equations with variable coefficients, *Chaos Solitons Fractals* 19, 2004, 847–851.

[18] T. Ozis, A. Yldrm. Applications of He,s semi-inverse method to the non-linear Schroinger equation. Computers and Mathematics and Applications. 54. 2007. 1039-1042.

[19] S. Guo, L. Mei and Z. Zhang. Time-fractional Gardner equation for ion-acoustic waves in negative-ion-beam plasma with negative ions and non-thermal non-extensive electrons. Physics of Plasma. 22. 2015.1-8.

[20] El-Wakil, S.A., Abulwafa, E.M., El-Shewy, E.K., Mahmoud, A.A.: Time-fractional KdV equation for plasma of two different temperature electrons and ion. Phys. Plasmas 18(9), 092116 (2011).

[21] El-Wakil, S.A., Abulwafa. Formulation and solution of space-time fractional Boussinesq equation. Nonlinear Dyn (2015) 80:167–175.

[22] Riewe, F.: Nonconservative Lagrangian and Hamiltonian mechanics. Phys. Rev. E **53**(2), 1890–1899 (1996).

[23] Riewe, F.: Mechanics with fractional derivatives. Phys. Rev. E **55**(3), 3581–3592 (1997).

[24] Agrawal, O.P.: Formulation of Euler–Lagrange equations for fractional variational problems. J. Math. Anal. Appl. **272**(1), 368–379 (2002).

[25] A. B. Malinowska, T. Odzijewicz, D. F.M. Torres , Advanced Methods in the Fractional Calculus of Variations. Springer. 2015.

[26] G. Jumarie. Lagragian mechanics of fractional order, Hamilton–Jacobi fractional PDE and Taylor's series non-differentiable functions. *Chaos, Solitons and Fractals*, **32**, 969–987 (2007).

[27] P.S. Bauer. Dissipative dynamical systems I. Proceedings of the National Academy of Sciences , 17:311–314, 1931.

[28] A. NNazari-Golshan and S. S. Nourazar. Effect of trapped electron on the dust ion acoustic waves in dusty plasma using time fractional modified Korteweg-de Vries equation. Physics of Plasmas.2013. 20(10). 103701.

[29] E. K. El-Shewy , A. A. Mahmoud , A. M. Tawfik , E. M. Abulwafa , and A. Elgarayhi . Space—time fractional KdV—Burgers equation for dust acoustic shock waves in dusty plasma with non-thermal ions. Chin. Phys. B 23, 070505 (2014)

[30] El-Sayed AMA, Gaber M. The Adomian decomposition method for solving partial differential equations of fractal order in finite domains. Phys. Lett. A. 2006; 359(20):175-182.

[31] S. Das, "Kindergarten of Fractional Calculus", (Lecture notes on fractional calculus-in use as limited prints at Dept. of Phys, Jadavpur University-JU; under publication at JU).

[32] J.-H.He, Homotopy perturbation technique, Comput Methods. Appl. Mech. Eng, 178. 1999. 257.

[33] H. Jafari and S. Momani, Solving fractional diffusion and wave equations by modified homotopy perturbation method, Phys. Lett. 2007 **A 370**. 388–396.

[34] A. Arikoglu and I. Ozkol. Solution of fractional differential equations by using differential transform method. Chaos, Soliton & Fractals. 2007. 54(5). 1473-1481.

[35] B. Lu . Bäcklund transformation of fractional Riccati equation and its applications to nonlinear fractional partial differential equations. Physics Letters A. 376. (2012). 2045-2048.

[36] Wu, G.-C., Lee, E.W.M.: Fractional variational iteration method and its application. Phys. Lett. A **374**(25), 2506–2509 (2010).

[37] Wu, G-c: A fractional variational iteration method for solving fractional nonlinear differential equations. Comput. Math. Appl. **61**(8), 2186–2190 (2011).

[38] S. Momani, Z. Odibat, Analytical solution of a time-fractional Navier–Stokes equation by Adomian decomposition method, Appl. Math. Comput. 177. (2006). 488–494.

[39] E. Yomba. The extended Fan's sub-equation method and its application to KdV–MKdV, BKK and variant Boussinesq equations. Physics Letters A 336 (2005) 463–476

[40] S. Zhang and H. Q. Zhang, Fractional sub-equation method and its applications to nonlinear fractional PDEs, Phys. Lett. A. 2011. 375. 1069.

[41] U. Ghosh, S. Sengupta, S. Sarkar, S. Das. Analytical solution with tanh-method and fractional sub-equation method for non-linear partial differential equations and corresponding fractional differential equation composed with Jumarie fractional derivative. Int. J. Appl. Math. Stat.; 2016. 54 (3). 11-31.

[42] U. Ghosh, S. Sarkar, S. Das. Analytical Solutions of Classical and Fractional KP-Burger Equation and Coupled KdV Equation. CMST. 2016. 22(3). 143-152.

[43] J. F. Alzaidy. The fractional sub–equation method and exact analytical solutions for some non-linear fractional PDEs. American Jounal of Mathematical Analysis. 2013.1(1) 14-19.

[44] S. Zhang, J-Lin Tong, W. Wang. A generalized ($G'/G$) for the mKdV equation with variable coefficients. Physics Letters A.2008 . 372. 2254-2257.

[45] Jumarie, G.: Modified Riemann–Liouville derivative and fractional Taylor series of non-differentiable functions further results. Comput. Math. Appl. **51**(9–10), 1367–1376 (2006).

[46] M. Caputo, "Linear models of dissipation whose *q* is almost frequency independent," *Geophysical Journal of the Royal Astronomical Society*, 1967. 13 (5). 529–539.

[47] Jumarie G, Tables of some basic fractional calculus formulae derived from modified Riemann-Liouville derivative for non-differentiable functions; Applied Mathematics Letters 22 (2009) 378-385.

[48] D. Baleanu, , Muslih, S.I.: Lagrangian formulation of classical fields with in Riemann–Liouville fractional derivatives. Phys. Scr. **72**(1), 119–123 (2005).

[49] D. Baleanu, , Muslih, S.I.. On fractional varioational principles. J. Sabaiter et al (eds), Advance in fractional calculus: Theoretical developments and Applications in Physics and Engineering. Springer. 2007. 115-126.

[50] Mittag-Leffler, G. M.. Sur la nouvelle fonction *Eα* (*x*), *C. R. Acad. Sci. Paris*, (Ser. II) 137, 554-558 (1903).

[51] Khater AH, Moussa MHM, Abdul-Aziz SF. Invariant variational principles and conservation laws for some nonlinear partial differential equations with constant coefficients–II. Chaos, Solitons & Fractals 2003;15:1–13.

[52] G. Jumarie. On the derivative chain-rules in fractional calculus via fractional difference and their application to systems modelling . Cent. Eur. J. Phys. 2013 . 11(6) . 617-63.

[54] U Ghosh , S Sarkar , S Das. Solution of System of Linear Fractional Differential Equations with Modified Derivative of Jumarie Type, American Journal of Mathematical Analysis, 2015, Vol. 3, No.3, pp72-84.

[55] M. Wang, X. Li and J. Zhang, The ($G'/G$)-expansion method and travelling wave solutions of nonlinear evolution equations in mathematical physics, Phys. Lett. A 372 (2008) 4. 417–423.

[56] U Ghosh, S Sarkar, S Das. Fractional Weierstrass Function by Application of Jumarie Fractional Trigonometric Functions and its Analysis, Advances in Pure Mathematics, 2015, 5, pp 717-732.